\documentclass[12pt]{article}

\usepackage{amsfonts}
\usepackage{amsfonts}  
\usepackage{amssymb}
\usepackage{amsthm}
\usepackage{amsmath}
\usepackage{color}

\title{H\"older regularity of the gradient for  solutions of  fully nonlinear equations with sublinear first order terms.} 

\author{}

\date{}

\catcode`@=11 \@addtoreset{equation}{section} \catcode`@=12

\newtheorem{theo}{Theorem}[section]
\newtheorem{prop}[theo]{Proposition}
\newtheorem{rema}[theo]{Remark}

\newtheorem{cor}[theo]{Corollary}
\newtheorem{lemme}[theo]{Lemma}

\def\R{\mathbb  R}
\def\grad{\nabla}

\def\bh#1{\|#1\|_{{\mathcal C}^{o,\gamma}}}

\DeclareMathOperator*{\osc}{osc}
\DeclareMathOperator*{\Lip}{Lip}

\begin{document}
\maketitle

\section{Introduction}
In this paper we shall establish some  regularity results of solutions of a class of fully nonlinear equations, 
with a first order term which is sub-linear. This paper is a natural continuation of \cite{BD10} and \cite{IS}. We now state the main result
and we will contextualize it later.

In the whole paper $\alpha$ will denote a number $> -1$, and $b$ will denote a function defined e.g. on $\R^N \times \R^N$ such that 
 
\begin{description}
 \item{[b1]} There exist  $\beta_o, \beta_\infty \in (0, 1+\alpha)$, $\gamma \geq \sup (\beta_0, \beta_\infty)$  and $M_o>0$ such that $\bh{b(., q)}\leq M_o (|q|^{\beta_o}+|q|^{\beta_\infty})$.

\item{[b2]} There exist 
$Q$ and $C_b$ such that for all $p\neq 0$ and $|q|\leq Q|p|$ 
$$|b(x, p+ q)-b(x, p)|\leq C_b\sup(|p|^{\beta_\infty-1}, |p|^{\beta_o-1})|q|.$$
\end{description} 

 In the whole paper $\beta$ denotes  $\sup (\beta_o, \beta_\infty)$, and 
 $\|b\| =  (M_o+ C_b)$, choosing $M_o$ and $C_b$ minimal in the above conditions.

We shall consider the following equation
\begin{equation}\label{08}
|\nabla u|^\alpha (F(D^2 u)+ h(x) \cdot \nabla u) + b(x,\nabla u)= f(x).
\end{equation} 
Observe that  hypothesis [b1] includes the case where $b(x, \cdot)$ is not Lipschitz continuous near zero.

\noindent We have in mind the following {\bf Examples:} 

1)$b(x,q)=b(x)|q|^\beta$ with $\beta\in (0,1+\alpha)$ and $b(x)$ H\"older continuous,

2) $b(x,p)=b(x)\cdot p |p|^{\beta-1}$ with $\beta\in (0,1+\alpha)$ and $b(x)$ H\"older continuous,

\noindent more in general

3)   $b(x,q)=\sum_{i=1}^k b_i(x) |q|^{\beta_i}+\sum_{j=1}^l b_j(x)\cdot q |q|^{\beta_j-1}$
  where $b_i$ are H\"older continuous functions and $\beta_i\in (0, 1+\alpha)$. 

$B_r$ indicate a ball of radius $r$, without loss of generality we will suppose that it is centered at the origin.
The main result is the following:
\begin{theo}\label{th1}
Let $\alpha>-1$. Suppose that 
$F$ is uniformly elliptic, that $h$ is a Lipschitz continuous function  which, when $\alpha >0$,  satisfies 
\begin{equation}\label{eqh2}(h(x)-h(y))\cdot (x-y)\leq 0
\end{equation}  
and $b$ satisfies [b1] and [b2]. Let $f\in{\cal C}(B_1)$, and if $\alpha <0$, $f$ is $|\alpha|$-H\"older continuous.  For any $u$,  
bounded viscosity solution  of  (\ref{08}) in $B_1$ and for 
 any $r<1$, there exist $\gamma^\star = \gamma^\star (\|f\|_\infty, \Lip(h),\beta_o, \beta_\infty,\gamma)\in(0,1)$  and 
$C = C(\gamma^\star)$ such that 
    $$\|u\|_{{\cal C}^{1, \gamma^\star} (B_r)} \leq  C\left(\|u\|_\infty+ \|b\|^{1\over 1+ \alpha-\beta_o}+\|b\|^{1\over 1+ \alpha-\beta_\infty}+\|f\|^{1\over 1+\alpha} _\infty\right).$$
    \end{theo}

Recall that Imbert and Silvestre in \cite{IS} proved an interior H\"older regularity for the gradient of the solutions
of
$$
|\grad u|^\alpha F(D^2 u) =f(x)
$$
when $\alpha\geq 0$. Their proof relies on a priori Lipschitz bounds, rescaling and an improvement of flatness Lemma and
on the classical regularity results of Caffarelli, and Evans  \cite{Ca,CaC} and \cite{E} for uniformly elliptic equations.

Following their breakthrough, in \cite{BD10}, we proved the same interior regularity for solutions of equation (\ref{08})
when $\alpha\geq 0$ and $b\equiv 0$. We also proved $C^{1,\gamma}$ regularity up to the boundary if the boundary datum
is sufficiently smooth. 
Recall that our main motivation to investigate the regularity of these solutions i.e. the simplicity of the first eigenvalue associated
to the Dirichlet problem for $|\nabla u|^\alpha (F(D^2 u)+ h(x) \cdot \nabla u)$, required continuity of the gradient up to the boundary.

When $\alpha\in (-1,0)$, in \cite{BD9} we proved ${\cal C}^{1, \gamma}$ regularity for solutions of the Dirichlet problem, 
using a fixed point argument which required global Dirichlet conditions on the whole boundary. 
So one of the question left open was: is the local regularity valid for $\alpha<0$?

Observe that the following holds:
\begin{prop}\label{prop0}
Suppose that for $\alpha \in (-1,0)$, $u$ is a viscosity solution of  (\ref{08}).
Then $u$ is a viscosity solution of 
\begin{equation}\label{24}F(D^2 u)+ h(x) \cdot \nabla u +|\grad u|^{-\alpha}b(x, \nabla u) =f(x)|\grad u|^{-\alpha}.
\end{equation}
\end{prop}
The proof  is postponed to the appendix. 
This means that, one way to answer the above question is to look at the regularity of solutions
of 
$$F(D^2 u)+ h(x) \cdot \nabla u +f(x)|\grad u|^{-\alpha}=0$$
and clearly $b(x,\grad u) :=f(x)|\grad u|^{-\alpha}$ when $f$ is H\"older continuous is just one of the example above.
So the case $\alpha<0$ reduces to the case $\alpha=0$ in (\ref{08}).
For completeness sake we have decided to treat a more general case i.e. (\ref{08}), but for clarity reasons
we shall first detail the proof in the case $\alpha=0$ (see section 3 below) and then explain the differences to bring when $\alpha\geq 0$ (section 4).
Proposition \ref{prop0} also explains why we treat sublinear terms that may have different growths in the gradient.

A consequence of Theorem \ref{th1} is the following further regularity of the solutions:
\begin{cor}
Under the assumptions in Theorem \ref{th1} and for $\alpha \leq 0$, assuming in addition that $f$ is H\"older continuous and 
$F$ is convex or concave the solutions 
of (\ref{08}) are ${\cal C}^2$.
\end{cor}
{\em Proof} By Proposition \ref{prop0} it is sufficient to consider $\alpha=0$.  We will check that if $u$ is solution of (\ref{08})  then it is a solution of 
$F(D^2v) =  -h(x) \cdot \nabla u-b(x, \nabla u)+ f(x)$. 
Let  $\varphi$  be test function at some point $\bar x$, so  that in the super solution case 
$(u-\varphi)(x) \geq (u-\varphi)(\bar x) = 0.$
But since $u$ is ${\cal C}^1$, it implies that
$\nabla \varphi (\bar x) = \nabla u(\bar x)$. Hence, by definition of viscosity super solution
\begin{eqnarray*}F(D^2 \varphi) (\bar x) &\leq&  -h(\bar x) \cdot \nabla \varphi( \bar x) -b(\bar x, \nabla \varphi (\bar x))+ f(\bar x) \\
&=& -h(\bar x) \cdot \nabla u(\bar x)-b(\bar x, \nabla u(\bar x))+ f(\bar x).
\end{eqnarray*}
Proceed similarly to prove that it is a sub solution. This ends the proof. 

\begin{rema} From the proof of this corollary it is clear that the solutions of (\ref{08}) are as smooth as the solutions of
$$F(D^2v)=g(x)$$
when $g$ is H\"older continuous, see e.g. \cite{NV}, \cite{Nir}, \cite{CaC1}.
\end{rema}

Even for $F(D^2u)=\Delta u$ it would be impossible to mention all the work that has been done on equation of the form
$$F(D^2u)+|\grad u|^p=f(x).$$
Interestingly most of the literature,  is concerned with the case $p>1$. In particular the  so called natural growth i.e. $p=2$ has been much studied in variational contexts and 
the behaviors are quite different when $p\geq2$ or $1<p<2$.  We will just mention the fundamental papers of Lasry and Lions \cite{LL} and
Trudinger \cite{Tams}.
And more recently the papers of Capuzzo Dolcetta, Leoni and Porretta \cite{clp} and Barles, Chasseigne, Imbert \cite{BCI}. In the latter
the H\"older regularity of the solution is proved for non local uniformly elliptic operators, and with lower order terms that may  be sublinear.

In the last section we briefly treat the boundary regularity so that the full result obtained in this paper is the following.
    \begin{theo}\label{th01}Suppose that $\Omega$ is  a bounded ${\cal C}^2$ domain of $\R^N$ and $\varphi\in {\cal C}^{1, \gamma_o} (\partial \Omega)$.
Under the assumptions of Theorem \ref{th1}, 
for any $u$,  
bounded viscosity solution  of 
\begin{equation}\label{(1)} \left\{ \begin{array}{lc}
|\grad u|^\alpha\left( F(D^2 u) +h(x)\cdot\grad u\right)+b(x,\grad u)= f(x) & {\rm in} \ \Omega\\
    u = \varphi & {\rm on} \ \partial \Omega
    \end{array}\right.
\end{equation}
there exist $\gamma = \gamma (\|f\|_\infty,\Omega, \Lip(h),\|b\|, \gamma_o, \gamma,\beta_o, \beta_\infty)$  and 
$C = C(\gamma)$ such that 
    $$\|u\|_{{\cal C}^{1, \gamma} (\overline{\Omega})} \leq  C\left( \|\varphi\|_{ {\cal C}^{1, \gamma_o} (\partial \Omega)}+  \|u\|_\infty+\|b\|^{1\over 1+ \alpha-\beta_o}+\|b\|^{1\over 1+ \alpha-\beta_\infty}+\|f\|_\infty\right).$$
\end{theo}

\section{Preliminaries and comparison principle.}
Let $S^{N}$ denote the symmetric $N\times N $ matrices.
In the whole paper $F(D^2u)$ indicates a uniformly elliptic operator satisfying, for some 
$0<\lambda\leq\Lambda$
$$\lambda {\rm tr} N\leq F(M+N)-F(M)\leq \Lambda  {\rm tr} N$$
for any $M\in S^{N}$ and any $N\in S^{N}$ such that $N\geq 0$.
The constants appearing in the estimates below often depend on $\lambda$ and $\Lambda$, but we will not specify
them explicitly when it happens.

We begin with a comparison principle  in the case $\alpha = 0$, that will be needed in section 5.   
\begin{theo} \label{thcompa} We suppose  that  $b$  satisfies [b1], [b2] with $\alpha=0$, and $h$ is Lipschitz continuous.           
Suppose that $j$ is some increasing function.
Suppose that $\Omega$ is an open  bounded set in $\R^N$, and that 
$u$ and $v$ are respectively   sub- and supersolutions of 
             $$F(D^2 u) +h(x)\cdot \nabla u  +b(x,\nabla u) -j(u)\leq f(x)$$
             $$ F(D^2 v)+  h(x)\cdot \nabla u + b(x,\nabla v) -j(v)\geq g(x)$$
with $u\geq v$ on $\partial \Omega$.  Suppose that $f \leq  g$ in $\Omega$,  $f$ and $g$ being continuous. 
Then  $u\geq v$ in $\Omega$. The same conclusion holds when $j$ is non decreasing but $f< g$ in $\Omega$. 
\end{theo}

{\em Proof of Theorem \ref{thcompa}.}
We are in the hypothesis of the classical comparison principle for viscosity solutions (see \cite{Users})
since  
$$|(h(x)-h(y))\cdot a(x-y)|\leq \Lip(h) a|x-y|^2$$
and, as soon as $\gamma \geq \beta$, for $a$ large,
\begin{eqnarray*}
 |(b(x, a(x-y))-b(y, a(x-y)) |&\leq& C_b(|x-y|^{\beta_o+ \gamma} a^{\beta_o}+|x-y|^{\beta_{\infty}+ \gamma} a^{\beta_{\infty}})\\
 & =& \omega(a|x-y|^2).
 \end{eqnarray*}

\begin{rema}
It is clear that the comparison principle fails when $j(.)=c.te$ and $f\equiv g$. For example, when $\Omega$ is the ball of radius 1,
then $u\equiv 0$ and $u(x)=C(1-|x|^\gamma)$ with $\gamma=\frac{2-\beta}{1-\beta}$ and $C=\gamma^{-1}(\gamma+N-2)^{\frac{1}{\beta-1}}$ 
are both solutions of equation
$$ \left\{ \begin{array}{lc}
\Delta u +|\grad u|^\beta=0 & {\rm in} \ \Omega\\
u=0 & {\rm on} \ \partial\Omega.
\end{array}\right.
$$        
\end{rema}

\begin{rema}
In the case $\alpha \neq 0$, the conclusion of Theorem \ref{thcompa} still hold if one considers instead $|\grad u|^\alpha(F(D^2 u) +h(x)\cdot \nabla u) +b(x,\nabla u) -j(u)\leq f(x)$ and the proof follows the arguments used in \cite{BD1}.
\end{rema}

\section{Interior regularity results.}
This section is devoted to the case $\alpha=0$ i.e. we shall study the local regularity of solutions of
\begin{equation}\label{31}
F(D^2 u) +  h(x)\cdot\grad u+b(x, \nabla u) = f(x).
\end{equation}

\begin{theo}\label{int}
Suppose that $h$ is a Lipschitz continuous function, and $b$ satisfies [b1], [b2] with $\alpha=0$. 
Let $f\in{\cal C}(B_1)$.  For any $u$,  
bounded viscosity solution  of (\ref{31}) in $B_1$
and for any $r<1$ there exist
$$\gamma^\star= \gamma^\star (\|f\|_\infty,  \Lip(h),\beta_o,\beta_\infty)\ \mbox{and}\ 
C = C(\gamma^\star)$$ 
such that 
    $$\|u\|_{{\cal C}^{1, \gamma^\star} (B_r)} \leq  C\left(\|u\|_\infty+ \|b\| ^{1\over 1-\beta_o}+ \|b\|^{1\over 1-\beta_\infty}+\|f\|_\infty\right).$$
    \end{theo}

\subsection{Lipschitz regularity.}                         
We now state and prove a   uniform Lipschitz regularity result.
It will be needed in the proof of the H\"older continuity of the
gradient in section 3.2.
\begin{lemme}\label{lem1}
Suppose that $H:B_1\times \R^N\rightarrow \R$ is  such that

$H(.,0)$ is bounded in $B_1$ 
and there exist $\beta_1, \beta_2 \in (0,1]$ and $C>0$ such that for all $q\in \R^N$, 
$$|H(x,q)-H(x,0)|\leq C(|q|^{\beta_1}+|q|^{\beta_2}).$$ 
Then there exists $C_o$ such that if $C<C_o$, any  bounded solution $u$ of 
       
$$F(D^2 u) + h(x) \cdot \nabla u + H(x, \nabla u) = f(x)\ \mbox{in}\ B_1$$
is Lipschitz continuous in $B_r$, for $r<1$  with some  Lipschitz constant   
depending on $r$, $\|f\|_\infty$, $C_o$ and $\|H(. 0)\|_\infty$.
\end{lemme}
{\em Proof of Lemma \ref{lem1}}
The proof proceeds as in \cite{IS}, \cite{BD10}. We outline it here,
in order to indicate the changes that need to be done. 

Let $r<r'<1$ and $x_o\in B_r$,  we consider on $B_{r'}\times B_{r'}$ the function
$$\Phi (x, y) =u(x)-u(y)- L_1\omega(|x-y|)- L_2 |x-x_o|^2- L_2|y-x_o|^2$$
where $\omega(s)=s-w_os^{\frac{3}{2}}$ for $s\leq (2/3w_o)^2$ and continued constantly; 
here $w_o$ is chosen in order that $(2/3w_o)^2>1$.

The scope is to prove that, for $L_1$ and $L_2$ independent of $x_o$,  chosen large enough,  
\begin{equation}\label{phi}
\Phi(x, y) \leq 0\ \mbox{on}\ B_{r}^2.
\end{equation}
This will imply   that $u$ is  Lipschitz continuous on $B_r$ by taking 
$x= x_o$, and letting $x_o$ vary.

So we begin to choose $L_2 > {8\sup u\over (r^\prime-r)^2}$ and later we shall  suppose that $L_1>> L_2$. 

Suppose by contradiction that $\Phi(\bar x,\bar y)=\sup\Phi(x,y)>0$. By the hypothesis on $L_2, $  $(\bar x,\bar y)$ is 
in the interior of $B_{r}^2$.

Let $L_1\geq 12\Lip(h)\sup u$. If $|\bar x-\bar y|\geq \frac{1}{2\Lip(h)}$ then $\Phi(\bar x,\bar y)\leq 0$. So we can suppose that  
\begin{equation}\label{18}|\bar x-\bar y|\leq \frac{1}{2\Lip(h)}.
\end{equation}

Proceeding in the calculations as in \cite{BCI} (see also \cite{IS} and \cite{BD5})
we get that if (\ref{phi})  is not true then, using the standard definition of closed semi-jets $\overline{J}^{2,+}$, $\overline{J}^{2,-}$ , there exist 
$$(q_x, X) \in \overline{J}^{2,+} u(\bar x), \  (q_y, -Y) \in \overline{J}^{2,-} u(\bar y)$$
with $q_x = L_1 \omega^\prime (|x-y|) {x-y\over |x-y|} + 2L_2(x-x_o)$  and
 $q_y = L_1 \omega^\prime (|x-y|) {x-y\over| x-y|} - 2L_2(y-x_o)$.  
If $\frac{L_2}{L_1}=o(1)$, then there exists $\kappa_1>0$ and $\kappa_2>0$ 
such that
$${\cal M}^+(X+Y)\leq -\kappa_1 L_1$$
and $|q_x|, |q_y|\leq \kappa_2 L_1$.
Using (\ref{18}), 

\begin{eqnarray*}
f(\bar x)&\leq& h(\bar x)\cdot q_x+H(\bar x, q_x)+F(X)\\
&\leq& h(\bar x)\cdot q_x+H(\bar x, q_x)+F(-Y)+ {\cal M}^+(X+Y)\\
&\leq& f(\bar y) - \kappa_1 L_1 + \|h\|_\infty L_2 + \Lip(h) \kappa_2 L_1 |\bar x-\bar y|\\
&&+\|H(., 0)\|_\infty+C  \sum_{i=1,2}(|q_x|^{\beta_i}+|q_y|^{\beta_i}).
\end{eqnarray*}
The terms $\|H(., 0)\|_\infty$ and $\|h\|_\infty L_2$ are $o(L_1)$, while for $C_o\leq\frac{\kappa_1}{16\kappa_2}$ 
\begin{eqnarray*}\Lip(h) \kappa_2 L_1 |\bar x-\bar y|+C  \sum_{i=1,2}(|q_x|^{\beta_i}+|q_y|^{\beta_i})&\leq& \frac{\kappa_1 L_1}{2}+C_o(1+4\kappa_2L_1)\\
&\leq&
\frac{3\kappa_1 L_1}{4}+C_o.
\end{eqnarray*}

In conclusion we have obtained that
$f(\bar x)\leq f(\bar y) -\frac{\kappa_1 L_1}{4} +o(L_1)$.
This is  a contradiction for $L_1$ large.

\begin{cor}\label{conv3}
Suppose that $(f_n(x))_n$ and $(H_n(x,0))_n$ are sequences converging uniformly respectively to  $f_\infty(x)$ and $H_\infty(x)$ 
on any compact subset of $B_1$,  and that there exist $\beta_1$ and $\beta_2 $ both in $(0,1]$,such that for all $q\in \R^N$, 
\begin{equation}\label{mm}
 \left\vert H_n(x,q)-H_n(x,0)\right\vert\leq \epsilon_n (|q|^{\beta_1}+|q|^{\beta_2})
\end{equation}
with $\epsilon_n\rightarrow 0$. 
Let $u_n$  be a sequence of solutions of
$$F(D^2u_n)+h(x)\cdot\grad u_n+H_n(x,\grad u_n)=f_n(x)\ \mbox{in }\ B_1.$$
If $\|u_n\|_\infty$ is a bounded sequence, then up to subsequences,
$u_n$ converges, in any compact subset of $B_1$, to $u_\infty$ a solution of the limit equation
$$F(D^2u_\infty)+h(x)\cdot\grad u_\infty+H_\infty(x)=f_\infty(x)\ \mbox{in }\ B_1.$$
\end{cor}

\subsection{H\"older regularity of the gradient.}
We will follow the line of proof in \cite{IS} and \cite{BD10}, which consists of three steps.

\noindent {\bf Step 1. Rescaling} Let  
$$\kappa = \left(2\|u\|_\infty + \|f\|_\infty+ \left(\frac{\|b\|}{\epsilon_o} \right)^{1\over 1-\beta_o}+ \left(\frac{\|b\|}{\epsilon_o} \right)^{1\over 1-\beta_\infty}\right)^{-1},$$
then  $u_\kappa = \kappa u$ satisfies 
     
$$ F(D^2 u_\kappa)+h(x)\cdot \nabla  u_\kappa+b_\kappa(x,\grad u_k)=f_\kappa (x)$$
where $b_\kappa(x,q)=\kappa b(x,\kappa^{-1}q)$, $f_\kappa = \kappa f$ .
Clearly $|f_\kappa|  \leq \epsilon_o$ and
observe that  $b_\kappa$ satisfies

$$\bh{b_\kappa(., q)}\leq M_o (\kappa^{1-\beta_o}|q|^{\beta_o}+ \kappa^{1-\beta_\infty} |q|^{\beta_\infty})\leq\epsilon_o(|q|^{\beta_o}+|q|^{\beta_\infty}).$$
Furthermore   for $|p|\neq 0$ and $|q|\leq Q|p|$ 

\begin{eqnarray*}
|b_\kappa(x,p+q)-b_\kappa(x,p)|&\leq& \kappa |b(x,\kappa^{-1}(p+q))-b(x,\kappa^{-1}p)|\\
&\leq & |q| C_b\sup( |\kappa^{-1}p|^{\beta_\infty-1}, |\kappa^{-1}p|^{\beta_o-1})\\
&\leq& \epsilon_o |q|\sup (|p|^{\beta_o-1}, |p|^{\beta_\infty-1}).
\end{eqnarray*}
This imply that we need to prove regularity only when $f$ and $b$ are "small" and it ends this step.

{\bf Step 2. Improvement of flatness.}

\noindent We now state and prove the improvement of flatness Lemma :
\begin{lemme}\label{flat}
There exist $\epsilon_o\in (0,1)$ and $\rho\in (0,1)$ depending on 
$(\beta_o,\beta_\infty, \|b\|, \\ 
\bh{h}, \|f \|_\infty, N)$ 
such that for any $p\in \R^N$ and for any viscosity solution $u$  of 
 $$F(D^2 u) +h(x)\cdot(\grad u+p)+ b(x,\nabla u+p) = f(x)\ \mbox{in}\ B_1$$
such that $\osc_{B_1} u\leq 1$ and $\|b\|+\|f\|_{L^\infty} \leq \epsilon_o$, 
there exists $q^\star \in \R^N$ such that 
$$\osc_{B_\rho} (u-q^\star \cdot x) \leq {1\over 2} \rho. $$
\end{lemme}
{\em Proof of Lemma \ref{flat}.} We argue by contradiction and then suppose   that, for any $n\in \mathbb{N}$, there exist
$p_n\in\R^N$, $f_n$ and $b_n$ such that 
$$\|b_n\|+\|f_n\|_{\infty} \leq \frac{1}{n}$$
and there exists $u_n$ solution of
$$F(D^2 u_n) + h(x)\cdot(\grad u_n+p_n)+b_n(x,\nabla u_n+p_n)  = f_n(x)\ \mbox{in}\ B_1$$
with $\osc_{B_1}u_n\leq 1$  and such that,  for any $\rho\in(0,1)$ and any $q^\star\in\R^N$
 
$$\osc_{B_\rho} (u_n-q^\star \cdot x) \geq {1\over 2} \rho. 
$$
Observe that $u_n-u_n(0)$ satisfies the same equation as 
$u_n$, it has oscillation $1$ and it is bounded, we can  then suppose that the sequence $(u_n)$ is bounded.
Suppose first that 
$|p_n|$ is bounded, so it converges, up to subsequences.
Let $v_n=u_n+p_n\cdot x$, it is a solution of
$$F(D^2v_n)+h(x)\cdot \grad v_n + b_n(x,\grad v_n)=f_n(x).$$
We can apply  Corollary \ref{conv3} with $H_n(x,q)=b_n(x,q)$,  
since (\ref{mm}) holds with $\beta_1=\beta_o$ and $\beta_2=\beta_\infty$.

\noindent Hence $v_n$ converges uniformly to $v_\infty$, a solution
of the limit equation
$$F(D^2 v_\infty) +h(y)\cdot (\grad v_\infty)=0\ \mbox{in}\ B_1;$$
furthermore $v_\infty$ satisfies, for any $\rho\in(0,1)$ and any $q^\star\in\R^N$,
\begin{equation}\label{caf}\osc_{B_\rho} (v_\infty-q^\star  \cdot x) \geq {1\over 2} \rho. \end{equation}
But this contradicts the classical
${\mathcal C}^{1,\alpha}$ regularity results see Evans \cite{E}, Caffarelli \cite{Ca} and Capuzzo Dolcetta, Vitolo \cite{cv} .

We suppose now that $|p_n|$ goes to infinity and, up to subsequences, let $p_\infty$ be the limit of  
$\frac{p_n}{|p_n|}$. By Ascoli's convergence theorem, we know that 
$\lim \frac{h(x)\cdot p_n}{|p_n|} := h(x)\cdot p_\infty$
converges in its ${\cal C}^{0, \gamma}$ norm.
We begin by proving that $h(x)\cdot p_\infty=0$.

 Indeed,  dividing the equation by $|p_n|$, 
the functions $v_n= \frac{u_n}{|p_n|}$ satisfy the 
following  equation 

$$F(D^2 v_n) + h(x) \cdot(\nabla v_n+\frac{p_n}{|p_n|}) +|p_n|^{-1}  b_n(x,|p_n|\nabla v_n+p_n)= {f_n(x)\over |p_n|}.$$
Let
$$H_n(x,q)=h(x)\cdot\frac{p_n}{|p_n|}+|p_n|^{-1}  b_n(x,|p_n|q+p_n).$$
We want to check that (\ref{mm}) holds.  Indeed by [b2] if $|q|\leq Q$  it holds with $\beta_1 = 1$  and $\epsilon_n=\|b_n\||p_n|^{\beta-1}$ . On the other hand if $|q|>Q$, using [b1],
\begin{eqnarray*}
|H_n(x,q)-H_n(x,0)|&\leq& |p_n|^{-1} (| b_n(x,|p_n|q+p_n)|+|b_n(x,p_n)| )\\
&\leq&  \|b_n\||p_n|^{\beta-1}(Q^{\beta-1}+\frac{1}{Q})|q|.
\end{eqnarray*}
Again we are in the hypothesis of  Corollary \ref{conv3},
we get that  
$(v_n)_n$  converges to $0$ which is  a 
solution of the limit equation i.e.
 
 \begin{equation}\label{eqh} h(x)\cdot p_\infty=0.
 \end{equation}
We have used the fact that $|p_n|^{-1}  b_n(x,p_n)$ tends to zero using [b1].
 

There are two cases,  suppose first that 
$h(x)\cdot p_n+b_n(x,p_n)$ is  bounded in ${\cal C}^{0, \gamma}$ for some $\gamma$.
Let $H_n(x,q)=h(x)\cdot p_n+b_n(x,q+p_n)$,  so that, up to a subsequence,  $H_n(x,0)$ converges uniformly to some function  
$H_\infty(x)$ while
$(u_n)_n$ is  a uniformly bounded sequence of solutions of 

        $$F(D^2 u_n) + h(x) \cdot \nabla u_n + H_n(x, \nabla u_n) = f_n(x).$$
Since condition
[b2] implies (\ref{mm}) with $\beta_1=\beta_2=1$ and $\epsilon_n=\|b_n\| |p_n|^{\beta-1}$, we can apply Corollary \ref{conv3}.
Up to a subsequence, $u_n$ converges to $u_\infty$ which is a solution of 

 $$F(D^2 u_\infty) + h(x) \cdot \nabla u_\infty + H_\infty(x)=0.$$
 Furthermore $u_\infty$ satisfies (\ref{caf}) , for any $\rho\in(0,1)$ and any $q^\star \in\R^N$.
 As in the case $p_n$ bounded, this contradicts the classical
${\mathcal C}^{1,\alpha}$ regularity results cited above.

\medskip
We are left to treat the case where $h(x)\cdot p_n+b_n(x,p_n)$ is unbounded.

Let $a_n=\bh{h(x)\cdot p_n+b_n(x,p_n)}$.
Observe first that, using [b1] and (\ref{eqh}), ${a_n\over |p_n|}$  goes to zero. 
  
Define  $v_n ={u_n\over a_n}$, then  for $H_n(x,q)=\frac{1}{a_n}{\left(h(x)\cdot p_n + b_n (x, a_n q+p_n)\right)}$
it satisfies 
$$F(D^2 v_n) + h(x)\cdot \nabla v_n+ H_n(x,\grad v_n)=\frac{1}{a_n}f_n(x).$$

We are in the hypothesis of Corollary \ref{conv3}, since, for $n$ large enough $|a_nq|\leq Q |p_n|$, 
hence by [b2]
$$\left\vert H_n(x, q)- H_n(x,0)\right\vert\leq  \|b_n\|  |p_n|^{\beta-1}|q|.$$
Using the compactness of $(v_n)$ and passing to the limit one gets that 

$$F(0) + h(x)\cdot \nabla (0) +  H_\infty(x)=0.$$
This yields a contradiction, indeed $ H_\infty(x)$ being the uniform limit of ${H_n (x,0)}$  it is of norm 1 in ${\cal C}^{0, \gamma}$.

This ends the proof of Lemma \ref{flat}.

\bigskip
{\bf Step 3. Conclusion.}
It is well known that it is enough to prove the following iteration process in order to prove Theorem \ref{int}.

\begin{lemme}\label{lem3a}
Suppose that  $\rho $, $\epsilon_o \in [0,1]$ and $b$  are as in Lemma \ref{flat} 
and suppose that $u$ is a viscosity solution of 
 $$F(D^2 u) + h(x)\cdot\grad u+b(x,\nabla u)  = f(x)\ \mbox{in}\ B_1$$
 with, $\osc u\leq 1$,  $\|f\|_\infty \leq \epsilon_o$, $M_o+ C_b \leq \epsilon_o$  then, there exists ${\gamma^\star} \in (0,1)$, such that for all $k>1$, $k\in \mathbb N$ there exists $p_k\in \R^N$ such that
\begin{equation}\label{kk2}
\osc_{B_{r_k}} (u(x)-p_k\cdot x) \leq r_k^{1+{\gamma^\star}} 
\end{equation}
 where $r_k:={\rho}^{k}$.
 \end{lemme}
{\em Proof.}
The proof is by induction and rescaling. For $k=0$ just take $p_k=0$. 
Suppose now that, for a fixe $k$, (\ref{kk2}) holds with some $p_{k}$. Choose 
${\gamma^\star}\in(0,1)$ such that $\rho^{\gamma^\star}>\frac{1}{2}$.

Define the function $u_k(x)= r_k^{-1-{\gamma^\star}} \left(u(r_k x) -p_k\cdot (r_k x)\right).$
By the induction hypothesis, $p_k$ is such that $\osc_{B_1} u_k\leq 1$
and $u_k$ is a solution of
$$F(D^2 u_k) +h(r_kx)r_k\cdot(\grad u_k+p_kr_k^{-{\gamma^\star}})+r_k^{1-{\gamma^\star}}b(r_k x,r_k^{\gamma^\star}(\grad u_k+p_kr_k^{-{\gamma^\star}}))  = r_k^{-{\gamma^\star}+1}f(r_kx)$$
Denoting by $b_k(x,q)=r_k^{1-{\gamma^\star}}b(r_kx,r_k^{\gamma^\star} q)$,we need to prove that $\| b_k\| \leq \epsilon_o$. Notice that

$$|{b_k(x,q)}|\leq r_k^{1-{\gamma^\star}}\epsilon_o(|r_k q|^{\gamma^\star \beta_o} +|r_k q|^{\gamma^\star \beta_\infty }) \leq \epsilon_o$$
so [b1] is satisfied with  $M_o = \epsilon_o$. 
  Suppose  now that $p\neq 0$ and $|q|\leq Q|p|$ i.e. $r_k|q|\leq Qr_k|p|$ we get
\begin{eqnarray*}
|b_k(x,p+q)-b_k(x,p)|&=&r_k^{1-\gamma^\star} |(b(r_k x,r_k^{\gamma^\star}(p+q))-b(r_k x,r_k^{\gamma^\star}p)|\\
&\leq & \epsilon_or_k|q|\sup((r_k^{\gamma^\star}|p|)^{\beta_\infty-1},(r_k^{\gamma^\star}|p|)^{\beta_o-1})\\
&\leq & \epsilon_o|q|\sup(|p|^{\beta_\infty-1},|p|^{\beta_o-1}).
\end{eqnarray*}
So we are in the hypothesis of Lemma \ref{flat}, hence there exists $q_k$ such that

$$\osc_{B_\rho} (u_k(x)-q_k \cdot x) \leq \frac{1}{2}\rho. $$
So that, for $p_{k+1}=p_k+q_kr_k^{{\gamma^\star}+1}$
$$
\displaystyle\osc_{ B_{r_{k+1}}} \left(u(x)-p_{k+1}\cdot x \right)
   \leq \frac{\rho}{2} r_k^{1+{\gamma^\star}} \leq r_{k+1}^{1+{\gamma^\star}}.
$$
This ends the proof of Lemma \ref{lem3a}, and in the same time the proof of Theorem \ref{int}.

\section{Regularity for  singular or degenerate  elliptic equations.}
We now prove Theorem \ref{th1} for $\alpha\neq 0$.
Here we use the definition of viscosity  solutions that we introduced in \cite{BD1}. Note that in the case $\alpha \geq 0$ that definition is equivalent to the classical definition.  

\bigskip
As seen in the introduction, through Proposition \ref{prop0} the case 
 $\alpha<0$,  reduces to the case of the previous section i.e. $\alpha=0$ replacing $b$ with 
 $\tilde b(x,p)=|p|^{-\alpha} b(x,p)-f(x)|p|^{-\alpha}$.

Hence suppose that $\alpha\geq 0$ and consider
\begin{equation}\label{345}|\grad u|^\alpha\left( F(D^2u)+h(x)\cdot \nabla u\right)+ b(x,\grad u)=f(x)\end{equation}
with $b$ satisfying [b1] and [b2].

\bigskip
In order to prove Theorem \ref{th1} we will follow the scheme used in section 3.2, just emphasizing the differences 
without writing down all the details. Let us note that the passage to the limit requires some compactness lemma, both for $p$ large and for $p$ bounded. This is done in Lemma \ref{holdalpha} and \ref{lipalpha} that we announce now and whose proofs are postponed to the appendix.

  \begin{lemme}\label{holdalpha}
Suppose that  
$H(.,0)$ is bounded in $B_1$ 
and there exist $\beta_1, \beta_2 \in (0,1+\alpha]$ and $C>0$ such that, for any $q\in\R^N$,
$$|H(x,q)-H(x,0)|\leq C(|q|^{\beta_1}+|q|^{\beta_2}).$$ 
Suppose that  $f $ is continuous and bounded. Then  the solutions of (\ref{345}) are H\"older continuous. 
\end{lemme}
              
\begin{lemme}\label{lipalpha}
Let $f$ and $H$ be  as in Lemma \ref{holdalpha} with  $\beta_1, \beta_2 \in (0,1]$.
Let $e\in \R^N$ of norm one.
There exist $a_o$, $C_o$ small such that for all $a<a_o$ and 
$C\leq C_o$, any solution of 
            $$|e+ a \nabla u|^\alpha (F(D^2 u) +h(x)\cdot \nabla u) +  H(x, \nabla u)= f(x)$$
is Lipschitz continuous. \end{lemme}
Of course these two lemmata will imply the corollaries 
             
\begin{cor}\label{compalpha1}
Suppose that $(u_n)$ is a bounded  sequence of solutions of
          
           $$|\nabla u_n|^\alpha \left(F(D^2 u_n) + h(x)\cdot \nabla u_n\right)+ H_n(x, \nabla u_n)= f_n(x)$$
and suppose that $f_n$   and $H_n(\cdot,  0)$ are uniformly  convergent respectively to $f_\infty$ and $H_\infty$.  Suppose that 
there exist $\beta_1, \beta_2 \in (0,1+\alpha]$  such that for all $q\in\R^N$
            \begin{equation}\label{mm1}
            |H_n(x,q)-H_n(x,0)|\leq \epsilon_n(|q|^{\beta_1}+|q|^{\beta_2}),
\end{equation}
with $\epsilon_n\rightarrow 0$.  Then $u_n$ converges up to subsequence to $u_\infty$ a solution of 
            $$|\nabla u_\infty|^\alpha (F(D^2 u_\infty)+ h(x) \cdot \nabla u_\infty) + H_\infty(x)= f_\infty(x).$$
\end{cor}
and 
     \begin{cor}\label{compalpha2}
Suppose that $(u_n)$ is a bounded sequence of  solutions of
         $$|e_n+ a_n\nabla u_n|^\alpha \left(F(D^2 u_n) + h(x)\cdot \nabla u_n\right)+ H_n(x, \nabla u_n)= f_n(x)$$
where $|e_n|=1$. Suppose that $H_n$ satisfies (\ref{mm1}) with  $\beta_1, \beta_2 \in (0,1]$ and
with $\epsilon_n$  a sequence which goes to zero.  Suppose that  $f_n$  and $H_n(., 0)$ converge uniformly respectively to $f_\infty$ and to $H_\infty$,  that $a_n$ tends to zero. Then $u_n$ converges, up to subsequences, to $u_\infty$ which satisfies 
            $$F(D^2 u_\infty) + h(x) \cdot \nabla u_\infty+  H_\infty (x) =f_\infty(x).$$
\end{cor}
We now use these results in order  to prove Theorem \ref{th1} using the same steps as in section 3.2.

{\bf Step 1.} For the rescaling $\tilde u=\kappa u$, choose
$$\kappa = \left( 2\|u\|_\infty +\left( {\|f\|_\infty\over \epsilon_o}\right)^{1\over 1+\alpha}  + \left(\frac{\|b\|}{\epsilon_o} \right)^{1\over 1+\alpha-\beta_o}+ \left(\frac{\|b\|}{\epsilon_o} \right)^{1\over 1+\alpha -\beta_\infty}\right)^{-1}.$$
Then, arguing as in step 1 in section 3.2,  it is enough to prove the regularity result for $\|b\|$ and $\|f\|_\infty$ small enough.

{\bf Step 2.} The statement and the proof of the improvement of flatness Lemma are similar to those of Lemma \ref{flat}. 
Let us detail a few passages.

First case  : the sequence $(p_n)$ is bounded. Since
$H_n(x,q) = b_n (x, q)$ satisfies (\ref{mm1}) with  $\beta_1 = \beta_o$ and $\beta_2  = \beta_\infty$, 
by Corollary \ref{compalpha1}, 
$u_n+ p_n \cdot x$ tends to $v_\infty$ a solution of 
$$|\nabla v_\infty |^\alpha  \left(F(D^2 v_\infty ) + h(x) \cdot \nabla v_\infty \right)=0\ \mbox{i.e.}\ F(D^2 v_\infty ) + h(x) \cdot \nabla v_\infty=0,$$
as seen in \cite{IS} and \cite{BD10}.
Furthermore,  for any $q^\star \in \R^N$
$\osc_{B_\rho} (u_n-q^\star \cdot x) \geq {1\over 2} \rho.
$
Passing to the limit, we get that  
$$\osc_{B_\rho} (v_\infty-\bar p \cdot x-q^\star \cdot x) \geq {1\over 2} \rho. 
$$ where $\bar p$ is the limit of a subsequence of $p_n$.   As in section 3.2, this contradicts the classical
${\mathcal C}^{1,\alpha}$ regularity results of Caffarelli \cite{Ca} and Evans \cite{E}.

Second case: $(p_n)_n$ is unbounded.
Up to subsequences, let $p_\infty$ be the limit of  
$\frac{p_n}{|p_n|}$. We begin by proving that $h(x)\cdot p_\infty=0$.

\noindent Indeed,  dividing the equation by $|p_n|^{\alpha+1}$, 
the functions $v_n= \frac{u_n+p_n\cdot x}{|p_n|}$ satisfy the 
following  equation

$$|\grad v_n|^\alpha (F(D^2 v_n) + h(x) \nabla v_n) +{1\over |p_n|^{1+\alpha} } b_n (x, |p_n|\nabla v_n)= {f_n(x)\over |p_n|^{\alpha+1}}.$$
Let $H_n(x,q)=|p_n|^{-1-\alpha}  b_n(x,|p_n|q)$;
using [b1],
\begin{eqnarray*}
|H_n(x,q)-H_n(x,0)|&\leq& |p_n|^{-1-\alpha} | b_n(x,|p_n|q)|\\
&\leq& \|b_n\| |p_n|^{\beta-1-\alpha}(|q|^{\beta_o}+|q|^{\beta_\infty}).
\end{eqnarray*}
Again we are in the hypothesis of  Corollary \ref{compalpha1},
hence $v_n$  converges to  $v_\infty(x)=p_\infty\cdot x$ a solution of the limit equation which becomes
 $h(x)\cdot p_\infty=0.$
 
 There are two cases to treat.  Suppose first that 
$h(x)\cdot p_n+b_n(x,p_n)|p_n|^{-\alpha}$ is  bounded in ${\cal C}^{0, \gamma}$ for some $\gamma$.
Up to a subsequence $h(x)\cdot p_n+ |p_n|^{-\alpha}b_n(x,p_n)$ converges uniformly to some function  $H_\infty(x)$.
Let 
$$H_n(x,q)=|e_n+ a_n q|^\alpha h(x)\cdot p_n+|p_n|^{-\alpha} b_n(x,q+p_n),$$
$u_n$ satisfies 
        $$|e_n+ a_n \nabla u_n|^\alpha (F(D^2 u_n) + h(x)\cdot  \nabla u_n) + H_n( x, \nabla u_n)= {f_n(x)\over |p_n|^\alpha}$$
Note that [b2] implies (\ref{mm1}) with $\beta_1 = 1= \beta_2$ and, for some $C_\alpha$  constant depending only  on $\alpha$, and
 $\epsilon_n=\|b_n\| |p_n|^{\beta-1-\alpha}+C_\alpha\sup_x{|h(x)\cdot p_n|\over |p_n|}.$
 
So, using Corollary \ref{compalpha2}, up to a subsequence, $u_n$ converges to $u_\infty$ which is a solution of 

 $$F(D^2 u_\infty) + h(x) \cdot \nabla u_\infty + H_\infty(x)=0.$$
 Furthermore $u_\infty$ satisfies (\ref{caf}) , for any $\rho\in(0,1)$ and any $q^\star \in\R^N$.
 As in section 3.2, this contradicts the classical
${\mathcal C}^{1,\alpha}$ regularity results of Caffarelli \cite{Ca} and Evans \cite{E}.

\medskip
We are left to treat the case where $h(x)\cdot p_n+|p_n|^{-\alpha}b_n(x,p_n)$ is unbounded.
Let $a_n=\bh{h(x)\cdot p_n+|p_n|^{-\alpha}b_n(x,p_n)}$. 
Observe first that,  using [b1] and   $h(x)\cdot p_\infty=0$, the sequence ${a_n\over |p_n|}$  goes to zero. 
  
Define  $v_n ={ u_n\over a_n}$, then  it satisfies 
$$|e_n+ {a_n\over |p_n|}\nabla v_n|^\alpha F(D^2 v_n) + h(x)\cdot \nabla v_n+ H_n(x,\grad v_n)={f_n(x)\over a_n|p_n|^\alpha}$$
with
$$H_n(x,q)=\frac{1}{a_n}{\left(|e_n+ {a_n\over |p_n|} q|^\alpha h(x)\cdot p_n +|p_n|^{-\alpha} b_n (x, a_n q+p_n)\right)}.$$
We are in the hypothesis of Corollary \ref{compalpha2}, since for $n$ large enough $|a_nq|\leq Q |p_n|$, 
by [b2]:

$$\left\vert H_n(x, q)- H_n(x,0)\right\vert\leq  \|b_n\|\left( |p_n|^{\beta-(1+\alpha)}+ C_\alpha\sup_x{|h(x)\cdot p_n|\over |p_n|}\right)  |q|. $$
Using the compactness of $(v_n)$ and passing to the limit one gets that 
$$F(0) + h\cdot \nabla (0) +H_\infty (x)=0.$$
This yields a contradiction: $ H_\infty(x)$ being the uniform limit of ${H_n (x,0)}$  it is of norm 1 in ${\cal C}^{0, \gamma}$. 
This ends the proof of the improvement of flatness Lemma.

{\bf Step 3. Conclusion}. The rescaling proceeds as in Lemma \ref{lem3a}. So this concludes the proof of Theorem \ref{th1}.

\section{Regularity results up to the boundary.}
In this section we shall prove that the solutions of  the equations  considered in the  previous sections 
are ${\cal C}^{1,\gamma}$ up to the boundary 
if the data are sufficiently regular. For the sake of clarity we will only  briefly treat  the  equations 
considered in  section 3, but it is not difficult to prove analogous results for the equations considered  
in section 4.

\bigskip
In the whole section, the function $a$ is in  ${\cal C}^2(\R^{N-1})$,  and, without loss of generality,    we suppose that
$a(0)=0$ and $\nabla a (0)=0$.  $B$ will denote an open bounded set in $\R^N$ that contains $(0,0)$; we will use the 
notations  $T$ for the portion of the boundary  $T = B \cap \{ y_N = a(y^\prime)\}$ and  $B^T= B \cap \{ y_N > a(y^\prime)\}$.

In order to prove the H\"older regularity of the gradient
for solutions of
\begin{equation}\label{eqqq}
\left\{\begin{array}{lc}
      F(D^2 u)+h(y)\cdot\grad u+ b(y,\nabla u) = f(y)& {\rm in} \  B^T,\\
         u=\varphi &\ {\rm on} \  \  T,
       \end{array}\right.
\end{equation}
we will use the same scheme as described before,  in the presence of the boundary
as it was done in \cite{BD10}.

We start with some Lipschitz regularity results up to the boundary.
\begin{lemme}\label{lemboundary} Let $\varphi\in {\cal C}^{1, \beta_o}(T) $.
Let $d$ be the distance to the hypersurface 
$\{ y_N = a(y^\prime)\}$. 
Then,  for all $r<1$, there exists   
$\delta_o$ depending on ($\|f\|_\infty$,$\beta_o$, $\beta_\infty$,  $\|h\|_\infty $, $\|b\|$, $r$,    $\Lip(\varphi)$),   
such that for all $\delta <\delta_o$   ,
if $u$ is a solution  of (\ref{eqqq})  such that $\osc u \leq 1$  
then it satisfies 
\begin{equation}\label{bnm} |u(y^\prime, y_N)-\varphi(y^\prime)| \leq {6\over \delta}  { d(y)}\ \mbox{ in}  
\ B_r(0)\cap\{y_N >a(y^\prime)\} .
\end{equation}
 \end{lemme}
{\em Proof.} 
We begin by choosing  $\delta< \delta_1$,  such that on $d(y) < \delta_1$  the distance 
is ${\cal C}^2$ and  satisfies 
$|D^2 d|\leq C_1$.  We shall also later choose  $\delta$ smaller in function of 
$(\|f\|_\infty, \|h\|_\infty , \|b\|,   N)$.  
               
We construct $w$ 
a super solution of 
\begin{equation}\label{ww}
 {\cal M}_{\lambda,\Lambda}^+ (D^2 w)+h(y)\cdot\nabla w- \|b\|(\nabla w|^{\beta_o}+ |\nabla w|^{\beta_\infty})  <-\|f\|_\infty,\  \end{equation}
 in $B \cap \{ y_N > a(y^\prime),\ d(y)< \delta\}$,
such that $w\geq u$ on $\partial (B \cap \{ y_N > a(y^\prime),\ d(y)< \delta\})$, here  ${\cal M}_{\lambda,\Lambda}^+$ denotes the well known
Pucci operator.
The comparison principle  in Theorem \ref{thcompa}, will imply that this inequality holds in $B \cap \{ y_N > a(y^\prime),\ d(y)< \delta\}$ and this will be exactly 
the upper bound in (\ref{bnm}).

Let $\psi$ be a solution of
$$\left\{\begin{array}{lc}
       {\cal M}_{\lambda,\Lambda}^+ (D^2 \psi) + h(y) \cdot \nabla \psi = 0& {\rm in} \  B^T\\
         \psi =\varphi & {\rm on} \  T.
       \end{array}\right.$$
It is well known that $\psi$ is ${\cal C}^{1, \gamma_o} (B\cap \{ y_N \geq a(y^\prime)\})\cap {\cal C}^2 ( B\cap \{ y_N >a(y^\prime)\})$.
Furthermore, we can choose $\psi$ such that $\|\psi\|_\infty\leq \|\varphi\|_\infty\leq 1$, 
$\|\nabla \psi\|_\infty \leq c \|\nabla \varphi\|_\infty$, for some constant $c$ which depends 
on  $B^T$,
see \cite{CaC}.
 Then, standards computations (\cite{BD10}) give that, for $m<1$ we can  choose $\delta$  small enough, such that the function
 $$w(y) = \left\{ \begin{array}{lc}
            {2\over \delta} {d(y) \over 1+ d^m(y)}+ \psi(y)&  {\rm for} \  |y  |  < r\\
             {2\over \delta} {d(y)\over 1+ d^m(y)} + {1\over  (1-r)^3} (|y|- r)^3+ \psi (y)&\ {\rm for} \   | y | \geq r,
             \end{array}\right.$$
 satisfies (\ref{ww}).

For the lower bound in (\ref{bnm}), we replace $w$ by $2\psi-w$. 
This ends the proof of  Lemma \ref{lemboundary}.

\begin{prop}\label{lip} Let $\varphi$ be a Lipschitz continuous function.
Suppose that $u$ satisfies (\ref{eqqq}).

For all $r<1$,   $u$ is $Lipschitz $ continuous   on 
$B_r \cap \{ y_N > a(y^\prime)\}$,  with some Lipschitz  constant depending on $(r, \beta_o, \beta_\infty, a,  \|f\|_\infty ,\|h\|_\infty, \Lip(h),  \|b\|,  \Lip(\varphi))$. 
\end{prop}

These two results imply the following corollary about convergence of solutions: 
\begin{cor}\label{conv6}
Suppose that $(f_n)_n$  and $H_n(x,0)$ are two sequences converging uniformly on
$\overline{B^T}$ to respectively $f_\infty$ and $H_\infty(x)$,  that there exist some 
$\beta_1$ and $\beta_2 $ both in $(0,1]$, such that for all $q\in \R^N$, 
\begin{equation}\label{mm5}
 \left\vert H_n(x,q)-H_n(x,0)\right\vert\leq \epsilon_n (|q|^{\beta_1}+|q|^{\beta_2})
\end{equation}
with $\epsilon_n\rightarrow 0$. 
Let $u_n$  be a sequence of solutions of
$$
\left\{\begin{array}{lc}
      F(D^2 u_n)+h(y)\cdot\grad u_n+ H_n(y,\nabla u_n) = f_n(y)& {\rm in} \  B^T,\\
         u=\varphi &\ {\rm on} \  \  T,
       \end{array}\right.
$$
such that $\|u_n\|_\infty$ is bounded, then up to subsequences,
$u_n$ converges,  for all $r<1$ in $ B_r\cap \{ y_N \geq a(y^\prime)\}$, to $u_\infty$ a solution of the limit equation
$$\left\{\begin{array}{lc}
      F(D^2 u_\infty)+h(y)\cdot\grad u_\infty+ H_\infty(y) = f_\infty(y)& {\rm in} \  B^T,\\
         u=\varphi &\ {\rm on} \  \  T.
\end{array}\right.
$$

\end{cor}

\bigskip

{\bf Step 2} The improvement of flatness Lemma is the following.

\begin{lemme}\label{lem1a} 
There exist $\epsilon_o>0$ and $\rho$  which depend on 
$(\beta_o, \beta_\infty,  \|D^2 a\|_\infty, \|b\|, \\
\Lip(h),\|\varphi\|_{{\cal C}^{1,\gamma_o}}) $ such that for any $p\in \R^N$ and $u$ a viscosity solution of  
$$\left\{ \begin{array}{lc}
 F(D^2 u)+h(y)\cdot(\grad u+p)+ b(y, \nabla u+p) = f(y)& \ {\rm in}\   B^T\\
 u+p\cdot y =\varphi & {\rm on} \ T.
 \end{array}\right.$$
  Then for all $x\in B$ such that $B_1(x) \subset B$, 
  $\osc_{B_1(x) \cap \{y_N >a(y^\prime)\}}  u \leq 1$, and $\|f\|_{L^\infty(B_1(x)\cap \{y_N >a(y^\prime)\})}\leq \epsilon_o$,   $\|b\|\leq \epsilon_o$,   there exists $q_{x, \rho} \in \R^N$ such that 
  
$$\osc_{B_\rho(x)\cap \{y_N >a(y^\prime)\}}  (u(y)-q_{x,\rho}\cdot y) \leq {\rho\over 2}.$$
\end{lemme}
{\em Proof of  Lemma \ref{lem1a}.}
We argue by contradiction and suppose that for all $n$,  
there exist  $x_n\in \overline{B}$  and $p_n\in \R^N$,  
$\|f_n\|_{L^\infty (\Omega )}\leq {1\over n} $, $\|b_n\| \leq {1\over n}$,  and  $u_n$ with $\osc   (u_n) \leq 1$
a solution of 
\begin{equation}\label{hhn}
\left\{ \begin{array}{lc}
 F(D^2 u_n)+ h(y) \cdot (\nabla u_n+p_n) + b_n (y,\nabla u_n+p_n)  = f_n(y) & {\rm in}\  B^T\\
 u_n(y)+ p_n\cdot y =\varphi(y^\prime)  & {\rm on} \  T,
 \end{array}\right.
 \end{equation}
such that  for any $q^\star\in\R^N$
\begin{equation}\label{bb}
\osc_{B_\rho(x_n)} (u_n(y)-q^\star  \cdot  y) >\frac{\rho}{2}.
\end{equation}
 Extract from $(x_n)_n$ a subsequence which converges to $ x_\infty \in \overline{B}\cap \{ y_N \geq a(y^\prime)\}$. 
 
We denote in the sequel $B_\infty= \cap_{n\geq N_1}  B_1(x_n)$,  which contains $B_\rho(x_\infty)$ as soon as $N_1$ is large enough.  

Observe that $u_n-u_n(x_n)$ satisfies the same equation as 
$u_n$, it has oscillation $1$ and it is bounded, we can  then suppose that the sequence $(u_n)$ is bounded.

The boundedness of $\|u_n\|_\infty$,   $\|\varphi_n\|_\infty$,  together with the boundary condition implies that
\begin{equation}\label{bord}
|p_n\cdot (y^\prime, a(y^\prime))|\leq C. 
\end{equation}
 {\bf The case where $T$ is not straight. }
 In that case, Proposition 3.6 in \cite{BD10} gives that $(p_n) $ is bounded. Then we will get a contradiction using Corollary \ref{conv6},
by showing that $v_n = u_n+ p_n\cdot y$ converges to a solution of the limit equation which will be smooth by the classical result of Caffarelli \cite{Ca} and Winter
\cite{NW} 
while at the same time it will satisfy 
$$\osc_{B_\rho(x_\infty)\cap \{y_N >a(y^\prime)\}}  (v(y)-(p^\prime_\infty +q_{x,\rho})\cdot y) \geq  {\rho\over 2},$$
where $p^\prime _\infty$ is the limit of a subsequence of $p_n^\prime$.

 \noindent {\bf The case where $T$ is straight i.e. $T = \{ y_N = 0\}\cap B$.}
Let $p_n=p_n^\prime+p_n^Ne_N$ then (\ref{bord}) implies that $(p_n^\prime)_n$ is bounded. 
If $(p_n^N)_n$ is bounded just proceed as above. 
So we suppose that $p_n^N$ is unbounded.  

Note that using the arguments as in the interior case one gets that $h_N (y)= \lim {h(y)\cdot p_n\over |p_n^N|} =0$
for every $y$  in $B_1\cap \{y_N > a(y^\prime)\}$.  

We define $H_n (y, q) =   b_n(y, p_n^Ne_N+ q)$.
We  suppose first that  $H_n(. ,0)$  is  bounded in ${\cal C}^{0, \gamma}$. Then $v_n = u_n+ p^\prime_n \cdot y$ satisfies 
 $$\left\{ \begin{array}{lc}
 F(D^2 v_n)+h(y)\cdot\grad v_n+ H_n(y,\nabla v_n)  = f_n(y)& \ {\rm in}\   B^T\\
 v_n =\varphi & {\rm on} \ T.
 \end{array}\right.$$
Since we are in the hypothesis of  Corollary \ref{conv6}, up to a subsequence $(v_n)_n$ converges to $v_\infty$. So
we conclude as in the non straight case.
 
Finally if 
 $a_n = \bh{H_n(.,0)}\rightarrow +\infty $. 
  we define $v_n = {u_n+p_n^\prime \cdot y\over a_n}$ and proceed as in the last case in step 2 of section 3.
 
{\bf Step 3.} In this case the conclusion is the following
 \begin{lemme}\label{lem2b}
Suppose that  $\rho $  and $\epsilon_o \in [0,1]$ and $B$ are as in Lemma \ref{lem1a} 
 and suppose that $u$ is a viscosity solution of (\ref{eqqq}), 
 with, $\osc u\leq 1$ and $\|f\|_\infty \leq \epsilon_o$,  then, there exists $\gamma^\star \in (0,1)$, such that for all $k$ and for all $x\in B$ such that $B_1(x) \subset B$, there exists $p_{x,k}\in \R^N$ for which  
 $$\osc_{B_{r_k}(x)\cap \{ y_N > a(y^\prime)\} } (u(y)-p_{x,k}\cdot y) \leq r_k^{1+\gamma^\star} $$
 where $r_k={\rho}^{k}$.
 \end{lemme}
The proof proceeds as in step 3 of section 3.2, together with the arguments used in Lemma 3.4   in \cite{BD10}.

 \section{Appendix : Proofs of Proposition \ref{prop0},  Lemmata \ref{holdalpha} and \ref{lipalpha}}
 {\em Proof of Proposition \ref{prop0}.}
We assume that $u$ is a supersolution  of (\ref{08}) and we 
want to prove that it is a supersolution  of (\ref{24}). Without loss of generality we shall prove it at the point $0$.  
If $u$ is constant around $0$, $D^2 u=0$ and $Du = 0$, so the conclusion is immediate. We suppose that 
there exists $M\in S$ such that 
\begin{equation}\label{mnb}
u(x)\geq u(0)+\frac{1}{2}\langle M(x),(x)\rangle +o(|x|^2)\end{equation} 
Let us observe first that one can suppose that $M$ is invertible, since if it is not, it can be replaced by  $M_n= M-{1\over n} I $ which satisfies (\ref{mnb}) and tends to $M$.

Let $k> 2$ and $R>0$ such that
          $$\inf_{ x\in B(0, R)} \left({u(x)-{1\over 2} \langle Mx, x\rangle+ |x|^k } \right)= u(0)$$
where the infimum is strict. 
We choose $\delta< R$ such that $(2\delta)^{k-2} < {\inf _i |\lambda_i (M)|\over 2k }$. 
Let $\epsilon $ be such that 
$$\inf_{ |x|> \delta }  \left(u(x)-{1\over 2} \langle Mx, x\rangle+  |x|^k \right)= u(0)+ \epsilon $$
and let $\delta_2< \delta $ and such that $k(2\delta)^{k-1}  \delta_2 + |M|_\infty (\delta_2^2 + 2 \delta_2 \delta) < {\epsilon\over 4}$. 
Then, for $x$ such that $|x|< \delta_2$,  
           \begin{eqnarray*}
           \inf_{ |y|\leq \delta} \{ u(y)-{1\over 2} \langle M(y-x), y-x\rangle +  |y-x|^k\} &\leq&      \inf_{|y|\leq \delta} \{ u(y)-{1\over 2} \langle My, y\rangle +  |y|^k\}+ {\epsilon\over 4}\\
           &=& u(0) + \epsilon/4
           \end{eqnarray*}
and on the opposite 
             \begin{eqnarray*}
             \inf_{ |y|> \delta} \{ u(y)-{1\over 2} \langle M(y-x), y-x\rangle +  |y-x|^k\} &\geq&      \inf_{|y|> \delta} \{ u(y)-{1\over 2} \langle My, y\rangle +  |y|^k\}- {\epsilon\over 4}\\
             &> &u(0) + {3\epsilon\over 4}.
             \end{eqnarray*}
Since the function $u$ is supposed to be non locally constant, there exist $x_\delta$ and $y_\delta$ in $B(0, \delta_2)$ such that 
        $$ u(x_\delta) > u(y_\delta) -{1\over 2} \langle M (x_\delta-y_\delta), x_\delta-y_\delta\rangle +  |x_\delta-y_\delta|^k$$
         and then the infimum $\inf _{y, |y| \leq \delta}\{u(y) -{1\over 2} \langle M (x_\delta-y), x_\delta-y\rangle + 
  |x_\delta-y|^k\}$ is achieved on some point $z_\delta$ different from $x_ \delta$. This implies that the function 
$$\varphi(z):=u(z_\delta ) + {1\over 2} \langle M(x_\delta- z), x_\delta-z) - |x_\delta-z|^k+ {1\over 2} \langle M(x_\delta-z_\delta), x_\delta-z_\delta \rangle +  |x_\delta -z_\delta|^k$$ 
touches $u$ by below at the  point $z_\delta$. 
But 
$$\grad\varphi(z_\delta)=M(z_\delta -x_\delta ) -k  |x_\delta-z_\delta|^{k-2} (z_\delta -x_\delta)\neq 0,$$    
indeed, if it was,  $z_\delta -x_\delta$ would be an eigenvector for the eigenvalue 
$k |x_\delta-z_\delta |^{k-2} $ which is supposed   to be strictly less than any eigenvalue of $M$. 

Since $u$ is a supersolution of (\ref{08}), multiplying by $|\grad\varphi(z_\delta)|^{-\alpha}$, we get
\begin{eqnarray*}
&&F\left( M - {d^2\over dz^2} (|x_\delta-z|^k)(z_\delta) \right)+ h(z_\delta)\cdot \grad\varphi(z_\delta) 
+ b\left(z_\delta, \grad\varphi(z_\delta)\right)|\grad\varphi(z_\delta)|^{-\alpha}\\
&\leq& f(z_\delta) |\grad\varphi(z_\delta)|^{-\alpha}.
\end{eqnarray*}

By passing to the limit we obtain 
               $$F(M) \leq 0$$
 which is the desired conclusion. 
               
We would argue in the same manner for sub-solutions. 

\bigskip

{\em Proof of Lemma \ref{holdalpha}.}
 Let $0< r< r^\prime <1$. We define  in $B_{r^\prime}$ when $x_o\in B_r$

$$\phi(x,y) = u(x)-u(y)-M|x-y|^{\gamma} -L |x-x_o|^2-L |y-x_o|^2.$$
We assume that $L > {8\sup u\over (r^\prime-r)^2}$ . We shall prove that $\phi\leq 0$ in $B_{r^\prime}$ this will implies,  
by taking first $x_o$ , and then exchanging $x$ and $y$ that $u$ is H\"older of exponent $\gamma$.
Proceeding by contradiction, as in the proof of Lemma \ref{lem1}, there exists  $(\bar x, \bar y)\in B_{r^\prime}$, $(q_x,X)\in \overline{J}^{2,+}(u(\bar x))$, $(q_y,-Y)\in \overline{J}^{2,-}(u(\bar y))$ such that
 $$tr(X+Y)\leq-\kappa_1 M |\bar x-\bar y|^{\gamma-2},\ q_x=\gamma M|\bar x-\bar y|^{\gamma-2}(\bar x-\bar y) -2L (\bar x-x_o)$$
and similarly for $q_y$, see \cite{BD10}.
Using the equation, and letting $C$ be a constant that may change from line to line, we get
    
                
\begin{eqnarray*}
             f(\bar x) &\leq&   |q_x|^\alpha F(X)+ h(x) \cdot q_x|q_x|^\alpha + H(\bar x, q_x)\\
             & \leq& |q_y|^\alpha  F(-Y)  + \lambda |q_x|^\alpha tr(X+Y)+(|q_x|^\alpha-|q_y|^\alpha)F(-Y)+ h(\bar y) \cdot q_y+\\
             &+ & \|h\|_\infty |q_x-q_y| |q_x|^\alpha + \|H(\cdot,  0)\|_\infty +C (\sum_{i=1,2} |q_x|^{\beta_i}+ |q_y|^{\beta_i}) + H(\bar y, q_y) \\
             &\leq & f(\bar y)-\kappa_1M^{1+\alpha} |\bar x-\bar y|^{(\gamma-2)+(\gamma-1)\alpha} +cL (M |\bar x-\bar y|^{\gamma-1})^{\alpha-1} (\|Y\| +\|h\|_\infty) \\
&&+C( \sum_i(M |\bar x-\bar y|^{\gamma-1})^{\beta_i}+\|H(\cdot,  0)\|_\infty\\
&\leq & -\kappa_1 M^{1+\alpha}  |\bar x-\bar y|^{(\gamma-2) + (\gamma-1)\alpha}+f(\bar y)+o(M^{1+\alpha}  |\bar x-\bar y|^{(\gamma-2)+ (\gamma-1)\alpha} ).
\end{eqnarray*}
We have used that $\|Y\|\leq CM |\bar x-\bar y|^{\gamma-2}$ and chosen $L= o(M)$. 
This is a contradiction for $M$ large enough.             

{\em Proof of Lemma \ref{lipalpha}.} As in the proof of Lemma \ref{lem1}, we want to prove that
there exists $L_1$ and $L_2$ such that, in $B_{r^\prime}$,
$$\phi(x,y) := u(x)-u(y)-L_1\omega(|x-y|) -L_2 |x-x_o|^2-L_2 |y-x_o|^2\leq 0,$$
for any $x_o\in B_r$ and $r< r^\prime <1$.  Proceeding by contradiction we obtain $(\bar x, \bar y)$, $(q_x,X)$, $(q_y,-Y)$
as in the proof of Lemma \ref{lem1}.

Let $a$ such $a \kappa_2 L_1 < {1\over 2}$.  We observe that 
 the inequality
            $$|e+ aq_x|^\alpha (F(X) + h(x)\cdot q_x)+ H(x, q_x) \geq f(\bar x)$$
implies that 
$$F(X) + h(x)\cdot q_x\geq 2^\alpha(-\|f\|_\infty -\|H\|_\infty (\cdot,  0)- C \kappa_2L_1)$$
  since  we can assume that $L_1>1$ and then one can  use $\beta_1 = \beta_2 = 1$ in the assumption on $H$ of Lemma \ref{lipalpha}. 
              
In the same manner 
               
$$F(-Y) + h(y)\cdot q_y\geq 2^\alpha (\|f\|_\infty +|H|_\infty (\cdot,  0)+ C \kappa_2L_1).$$
We now observe that 
               $$F(-Y) + h(y) \cdot q_x\leq F(-Y) + h(y)\cdot q_y+\|h\|_\infty |q_x-q_y|.$$
   We obtain that 
                 \begin{eqnarray*}
                 2^\alpha(-\|f\|_\infty -\|H (\cdot, 0)\|_\infty- C\kappa_2L_1)&\leq& F(X) + h(x)\cdot q_x\\
                 &\leq & F(-Y) + h(x)\cdot q_x-\kappa_1 L_1 \\
                 &\leq& F(-Y) + h(y)\cdot q_y +c \|h\|_\infty L_2 -\kappa_1 L_1\\
                 &\leq & 2^\alpha (\|f\|_\infty +\|H(\cdot, 0)\|_\infty + C\kappa_2L_1) \\
                 &&+c \|h\|_\infty L_2 -\kappa_1 L_1.
                 \end{eqnarray*}
This is clearly false for $L_1$ large enough as soon as $C$ is chosen small enough and $L_2$ is small with respect to $L_1$.

    \end{document}